\documentclass [a4 paper,12pt]{article}
\usepackage{mathrsfs}

\usepackage[usenames]{color}
\usepackage{bbm,amsmath,amssymb,amsfonts,latexsym,amsthm,graphics,epsfig,graphicx}
\usepackage{times}
\usepackage{indentfirst,graphicx,epsfig}

\newtheorem{lem}{Lemma}
\newtheorem{thm}[lem]{Theorem}
\newtheorem{cor}[lem]{Corollary}
\newtheorem{con}{Conjecture}
\newtheorem{rem}{Remark}
\newtheorem{pro}{Proposition}

\title{Derivatives and Real Roots of Graph Polynomials}
\author{\small Xueliang Li and Yongtang Shi\\
\small Center for Combinatorics and  LPMC\\
\small Nankai University, Tianjin 300071, China\\
\small Email: lxl@nankai.edu.cn,
 shi@nankai.edu.cn}
\date{ }

\begin{document}

\maketitle

\begin{abstract}
Graph polynomials are polynomials assigned to graphs. Interestingly, they also arise
in many areas outside graph theory as well. Many properties of graph polynomials
have been widely studied. In this paper, we survey some results on the derivative
and real roots of graph polynomials, which have applications in chemistry, control
theory and computer science. Related to the derivatives of graph polynomials,
polynomial reconstruction of the matching polynomial is also introduced.
\\
[2mm] Keywords: graph polynomial; derivatives; real roots; polynomial reconstruction\\
[2mm] AMS Subject Classification (2010): 05C31, 05C90, 05C35, 05C50
\end{abstract}

\section{Introduction}

Many kinds of graph polynomials have been introduced and extensively
studied, such as characteristic polynomial, chromatic polynomial,
Tutte polynomial, matching polynomial, independence polynomial,
clique polynomial, etc.

Let $G$ be a simple graph with $n$ vertices and $m$ edges, whose
vertex set and edge set are $V(G)$ and $E(G)$, respectively. The
{\it complement} $\overline{G}$ of $G$ is the simple graph whose
vertex set is $V(G)$ and whose edges are the pairs of nonadjacent
vertices of $G$. For terminology and notation not defined here, we
refer to \cite{BondyMurty}.

Denote by $A(G)$ the adjacency matrix of $G$. The {\it characteristic polynomial} of
$G$ is defined as
$$
\phi(G,x)=det(\lambda I-A(G))=\sum_{i=0}^{n}a_{i}^{n-i}.
$$
The roots $\lambda_{1}, \lambda_{2}, \ldots, \lambda_{n}$ of
$\phi(G,x)=0$ are called the {\it eigenvalues} of $G$. For more
results on $\phi(G,x)$, we refer to \cite{CDGT,CDS}.

Denote by $m(G,k)$ the {\it $k$-th matching number} of $G$ for $k \geq 0$. We assume
that $m(G,0)=1$. For $k\geq 1$, $m(G,k)$ is defined as the number of
ways in which $k$ pairwise independent edges can be selected in $G$. The {\it
matching polynomial} is defined as
$$\alpha(G,x)=\sum_{k\geq 0}(-1)^km(G,k)x^{n-2k}.$$
There is also an auxiliary polynomial $\alpha(G,x,y)$, which is
defined as
$$\alpha(G,x,y)=\sum_{k\geq 0}(-1)^km(G,k)x^{n-2k}y^k.$$
Note that
$\alpha(G,x,y)=y^{n/2}\alpha(G,xy^{-1/2})$. In view of this fact, we
may define an auxiliary polynomial of $\phi(G,x,y)$:
$$\phi(G,x,y)=y^{n/2}\phi(G,xy^{-1/2})=\sum_{k\geq 0}a_kx^{n-k}y^{k/2}.$$
Note that $\phi(G,x,y)$ is a polynomial in $y$ if and only if $G$ is
bipartite.

Denote by $n(G,k)$ the {\it $k$-th independence number} of $G$ for
$k \geq 0$. We assume that $n(G,0)=1$. For $k\geq 1$,
$n(G,k)$ is defined as the number of ways in which $k$ pairwise
independent vertices can be selected in $G$. The {\it independence
polynomial} is defined as
$$\omega(G,x)=\sum_{k\geq 0}(-1)^kn(G,k)x^{n-k},$$
which is also called independent set polynomial in \cite{HoedeLi} and
stable set polynomial in \cite{Stanley}. For more results on the
independence polynomials, we refer the surveys
\cite{LevitMadrescu,Trinks}.

Denote by $c(G,k)$ the $k$-th {\it clique number} of $G$ for $k \geq
0$. We assume that $c(G,0)=1$. For $k\geq 1$,
$c(G,k)$ is defined as the number of ways in which $k$ pairwise
adjacent vertices can be selected in $G$. Note that $c(G,1)=n$ and $c(G,2)=m$. The
{\it clique polynomial} is defined as
$$c(G,x)=\sum_{k\geq 0}(-1)^kc(G,k)x^{n-k}.$$
Note that the clique polynomial of a graph $G$ is exactly the
independence polynomial of the complement $\overline{G}$ of $G$,
i.e., $c(G,x)=\alpha(\overline{G},x)$. Obviously, we also have
$$c(G,x)+c(\overline{G},x)=\alpha(G,x)+\alpha(\overline{G},x).$$
The following results are easily obtained.
\begin{thm}
Let $G_1$ and $G_1$ be two vertex-disjoint graphs. Then we have
$$c(G_1\cup G_2,x)=c(G_1,x)+c(G_2,x)-1,\qquad \alpha(G_1\cup G_2,x)=\alpha(G_1,x)\cdot \alpha(G_2,x);$$
$$c(G_1+ G_2,x)=c(G_1,x)\cdot c(G_2,x),\qquad \alpha(G_1+ G_2,x)=\alpha(G_1,x)+ \alpha(G_2,x)-1.$$
\end{thm}

In \cite{HoedeLi}, the authors obtained the following similar
result.
\begin{thm}
Let $G_1$ and $G_2$ be two vertex-disjoint graphs with $n_1$ and
$n_2$ vertices, respectively. Then
$$c(G_1\times G_2, x)=n_2\cdot c(G_1,x)+n_1\cdot c(G_2,x)-(n_1+n_2+n_1n_2x)+1.$$
\end{thm}
For more properties on $c(G,x)$ and $\alpha(G,x)$, we refer to
\cite{HoedeLi}.

Many properties of graph polynomials have been widely studied. In this
paper, we survey some results on the derivative and real roots of graph
polynomials, which have applications in chemistry, control theory
and computer science. Related to the derivatives of graph
polynomials, polynomial reconstruction of the matching polynomial is
also introduced.

\section{Derivatives of graph polynomials}

The derivatives of the characteristic polynomial were examined by
Clarke \cite{Clarke} and the following result was showed.
\begin{thm}
Let $G$ be a simple graphs and $\phi(G,x)$ be the characteristic
polynomial of $G$. Then
$$\mathrm{d} \phi(G,x)/ \mathrm{d}x=\sum_{v\in V(G)}\phi(G-v,x).$$
\end{thm}
Gutman and Hosoya \cite{GutmanHosoya} got a similar result for the
matching polynomial.
\begin{thm}
Let $G$ be a simple graphs and $\alpha(G,x)$ be the matching
polynomial of $G$. Then
$$\mathrm{d} \alpha(G,x)/ \mathrm{d}x=\sum_{v\in V(G)}\alpha(G-v,x).$$
\end{thm}
One can get the first derivative of the independence polynomial and
clique polynomial, which have similar expressions as the
matching polynomial and characteristic polynomial. That is,
$$\mathrm{d} \omega(G,x)/ \mathrm{d}x=\sum_{v\in V(G)}\omega(G-v,x),\ \ \mathrm{d}
c(G,x)/ \mathrm{d}x=\sum_{v\in V(G)}c(G-v,x).$$
Although the four
first derivatives obey fully analogous expressions, their proofs
existing in the literatures, are quite dissimilar. Li and Gutman
\cite{LiGutman1995} provided a unified approach to all these
formulas by introducing a general graph polynomial.

Let $f$ be a complex-valued function defined on the set
of graphs $\mathcal{G}$ such that $G_1\cong G_2$ implies
$f(G_1)=f(G_2)$. Let $G$ be a graph on $n$ vertices and $S(G)$ be
the set of all subgraphs of $G$. Define $$S_k(G)=\{H: H\in S(G)\ and
\ |V(H)|=k\}, \ \ \ p(G,k)=\sum_{H\in S_k(G)}f(H).$$
Then, the general graph polynomial of $G$ is defined as
$$P(G,x)=\sum_{k=0}^np(G,k)x^{n-k}.$$
Actually, let $$f(H)=\left\{
 \begin{array}{ll}
 (-1)^{|V(H)|/2} &{\text{if $H$ is $1$-regular};}\\
0 &{\text{otherwise.}}
 \end{array}
 \right.
 $$ Then the resulting polynomial is the matching polynomial. Let $$f(H)=\left\{
 \begin{array}{ll}
 (-1)^{|V(H)|} &{\text{if $H$ is no edges};}\\
0 &{\text{otherwise.}}
 \end{array}
 \right.
 $$ Then the resulting polynomial is the independence polynomial. Let $$f(H)=\left\{
 \begin{array}{ll}
 (-1)^{r(H)}\cdot 2^{c(H)} &{\text{if all components of $H$ are 1- or 2-regular};}\\
0 &{\text{otherwise,}}
 \end{array}
 \right.
 $$ where $r(H)$ is the number of components in $H$ and $c(H)$ is the number of cycles in $H$.
 Then the resulting polynomial is the characteristic polynomial. Let  $$f(H)=\left\{
 \begin{array}{ll}
 (-1)^{|V(H)|} &{\text{if $H$ is a complete graph};}\\
0 &{\text{otherwise.}}
 \end{array}
 \right.
 $$ Then the resulting polynomial is the clique polynomial.

The following theorem was obtained by Li and Gutman in
\cite{LiGutman1995}.
\begin{thm}
For the graph polynomial $P(G,x)$ of $G$, we have
\begin{align*}
\frac{\mathrm{d}}{\mathrm{d}x} (P(G,x)) =\sum_{v\in
V(G)}P(G-v,x).
\end{align*}
\end{thm}

Furthermore, Gutman \cite{gutman1992, Gutman2} got the first derivative
formula for $\alpha(G,x,y)$:
$$\partial \alpha(G,x,y)/ \partial y=-\sum_{uv\in E(G)}\alpha(G-u-v,x,y).$$
To find an expression for $\partial \phi(G,x,y)/ \partial y$
of a bipartite graph was posed by Gutman as a problem in
\cite{Gutman0}. A solution of this problem was offered by Li and
Zhang \cite{LiZhang}.
\begin{thm}
For a bipartite graph $G$,
\begin{align}\label{eq1}
\partial \phi(G,x,y)/
\partial y=-\sum_{uv\in E(G)}\phi(G-u-v,x,y)-\sum_{C\subseteq G}
n(C)y^{n(C)/2-1}\phi(G-C,x,y),
\end{align}
where $C$ is a cycle, possessing $n(C)$ vertices.
\end{thm}

The above theorem was proved by using Sachs Theorem for the
coefficients of the characteristic polynomial and by verifying the
equality of the respective coefficients of the polynomials on the
left- and right-hand sides of Eq. \eqref{eq1}. In
\cite{GutmanLiZhang}, the authors put forward another route to Eq.
\eqref{eq1}, from which it become evident that Eq. \eqref{eq1} holds
for an arbitrary graph.

Moreover, if we define
$$P(G,x,y)=\sum_{i+j=n}p(G,k)x^iy^j,$$
then we can obtain
$$\frac{\partial P(G,x,y)}{\partial y}=ny^{-1}P(G,x,y)-xy^{-1}\sum_{v\in V(G)}P(G-v,x,y).$$

Derivatives of other graph polynomials have also been studied, such as the
cube polynomial \cite{BresarKlavzarSkrekovski}, the Tutte polynomial
\cite{Ellis-MonaghanMerino}, the Wiener polynomial
\cite{KonstantinovaDiudea}, etc.

\section{Polynomial reconstruction of the matching polynomial}

The derivative of a graph polynomial is related the problem of
polynomial reconstruction. This section aims to prove that graphs
with pendant edges are polynomial reconstructible and, on the other
hand, to display some evidence that arbitrary graphs are not, which
is given in \cite{LiShiTrinks}.

The famous (and still unsolved) reconstruction conjecture of Kelly
\cite{kelly1957} and Ulam \cite{ulam1960} states that every graph
$G$ with at least three vertices can be reconstructed from (the
isomorphism classes of) its vertex-deleted subgraphs.

With respect to a graph polynomial $P(G)$, this question may be
adapted as follows: Can $P(G)$ of a graph $G = (V, E)$ be
reconstructed from the graph polynomials of the vertex
deleted-subgraphs, that is from the collection $P(G_{-v})$ for $v
\in V$ ? Here, this problem is considered for the matching polynomial
of a graph. For results about the polynomial reconstruction of other
graph polynomials, see the article by Bre\v{s}ar, Imrich, and
Klav\v{z}ar \cite[Section 1]{bresar2005} and the references therein.
For additional results, see \cite[Section 7]{tittmann2011}
\cite[Subsection 4.7.3]{trinks2012c}.

The matching polynomial we considered here is the generating
function of the number of its matchings with respect to their
cardinality, denoted by $M(G, x, y)$, which is different from the
above $\alpha(G, x)$ and $\alpha(G, x, y)$. Let $G = (V, E)$ be a
graph. A {\it matching} in $G$ is an edge subset $A \subseteq E$,
such that no two edges in $A$ have a common vertex. The {\it matching
polynomial} $M(G, x, y)$ is defined as
$$
M(G, x, y) = \sum_{A \subseteq E \text{ is a matching}}
x^{\text{def}(G, A)} y^{|A|},
$$
where $\text{def}(G, A) = |V| - |\bigcup_{e \in A}{e}|$ is the
number of vertices not included in any of the edges of $A$. A
matching $A$ is a {\it perfect matching}, if its edges include all
vertices, that means if $def(G, A) = 0$. A {\it near-perfect
matching} $A$ is a matching that includes all vertices except one,
that means $def(G, A) = 1$. For more information about matchings and
the matching polynomial, see \cite{farrell1979b, gutman1977,
lovasz1986}.

For a graph $G = (V, E)$ with a vertex $v \in V$, $G_{-v}$ is the
graph arising from the {\it deletion} of $v$, i.e., arising by the
removal of $v$ and all the edges incident with $v$. The multiset of
(the isomorphism classes of) the vertex-deleted subgraphs $G_{-v}$
for $v \in V$ is the {\it deck} of $G$. The {\it polynomial deck}
$\mathcal{D}_P(G)$ with respect to a graph polynomial $P(G)$ is the
multiset of $P(G_{-v})$ for $v \in V$. A graph polynomial $P(G)$ is
{\it polynomial reconstructible}, if $P(G)$ can be determined from
$\mathcal{D}_P(G)$.

By arguments analogous to those used in Kelly's Lemma
\cite{kelly1957}, the derivative of the matching polynomial of a
graph $G = (V, E)$ equals the sum of the polynomials in the
corresponding polynomial deck.

\begin{pro}[Lemma 1 in \cite{farrell1987}]
Let $G = (V, E)$ be a graph. The matching polynomial $M(G, x, y)$
satisfies
\begin{align}
\frac{\delta}{\delta x} M(G, x, y) = \sum_{v \in V}{M(G_{-v}, x,
y)}.
\end{align}
\end{pro}

In other words, all coefficients of the matching polynomial except
the one corresponding to the number of perfect matchings can be
determined from the polynomial deck and thus also from the deck:
\begin{align}
m_{i, j}(G) = \frac{1}{i} \sum_{v \in V}{m_{i, j}(G_{-v})} \qquad
\forall i \geq 1,
\end{align}
where $m_{i, j}(G)$ is the coefficient of the monomial $x^i y^j$ in
$M(G,x,y)$.

Consequently, the (polynomial) reconstruction of the matching
polynomial reduces to the determination of the number of perfect
matchings.

\begin{pro} \label{prop:polynomial_reconstruction}
The matching polynomial $M(G, x, y)$ of a graph $G$ can be
determined from its polynomial deck $\mathcal{D}_M(G)$ and its
number of perfect matchings. In particular, the matching polynomial
$M(G, x, y)$ of a graph with an odd number of vertices is polynomial
reconstructible.
\end{pro}

Tutte \cite[Statement 6.9]{tutte1979} showed that the number of
perfect matchings of a simple graph can be determined from its deck
of vertex-deleted subgraphs and therefore gave an affirmative answer
on the reconstruction problem for the matching polynomial.

The matching polynomial of a simple graph can also be reconstructed
from the deck of edge-extracted and edge-deleted subgraphs
\cite[Theorem 4 and 6]{farrell1987} and from the polynomial deck of
the edge-extracted graphs \cite[Corollary 2.3]{gutman1992}. For a
simple graph $G$ on $n$ vertices, the matching polynomial is
reconstructible from the collection of induced subgraphs of $G$ with
$\lfloor{\frac{n}{2}}\rfloor + 1$ vertices \cite[Theorem
4.1]{godsil1981b}.

The following result is from \cite{LiShiTrinks} for simple graphs with pendant edges.

\begin{thm} \label{theo:pendant_perfect_matching}
Let $G = (V, E)$ be a simple graph with a vertex of degree $1$. Then, $G$
has a perfect matching if and only if each vertex-deleted subgraph
$G_{-v}$ for $v \in V$ has a near-perfect matching.
\end{thm}

As proved recently by Huang and Lih \cite{huang2014}, this statement
can be generalized to arbitrary simple graphs.

\begin{cor} \label{coro:forest_perfect_matching}
Let $G = (V, E)$ be a forest. Then $G$ has a perfect matching if and only
if each vertex-deleted subgraph $G_{-v}$ for $v \in V$ has a
near-perfect matching.
\end{cor}

Forests have either none or one perfect matching, because every
pendant edge must be in a perfect matching (in order to cover the
vertices of degree $1$) and the same holds recursively for the
subforest arising by deleting all the vertices of the pendant edges.
Therefore, from Proposition \ref{prop:polynomial_reconstruction} and
Corollary \ref{coro:forest_perfect_matching} the polynomial
reconstructibility of the matching polynomial follows.

\begin{cor}
The matching polynomial $M(G, x, y)$ of a forest is polynomial
reconstructible.
\end{cor}

On the other hand, arbitrary graphs with pendant edges can have more
than one perfect matching. However, Corollary
\ref{coro:forest_perfect_matching} can be extended to obtain the
number of perfect matchings. For a graph $G = (V, E)$, the number of
perfect matchings and of near-perfect matchings of $G$ is denoted by
$p(G)$ and $np(G)$, respectively.

\begin{thm} \label{theo:pendant_number_perfect_matching}
Let $G = (V, E)$ be a simple graph with a pendant edge $e = \{u,
w\}$ where $w$ is a vertex of degree $1$. Then we have
\begin{align}
&p(G) = np(G_{-u}) \leq np(G_{-v}) \qquad \forall v \in V \text{ and particularly} \\
&p(G) = \min{\{np(G_{-v}) \mid v \in V\}}.
\end{align}
\end{thm}

By applying this theorem, the number of perfect matchings of a
simple graph with pendant edges can be determined from its
polynomial deck and the following result is obtained as a corollary.

\begin{cor}
The matching polynomial $M(G, x, y)$ of a simple graph with a
pendant edge is polynomial reconstructible.
\end{cor}

While it is true that the matching polynomials of graphs with an odd
number of vertices or with a pendant edge are polynomial
reconstructible, it does not hold for arbitrary graphs.

There are graphs which have the same polynomial deck and yet their
matching polynomials are different. Although there are already
counterexamples with as little as six vertices, it seems that
nothing has been published before in connection with the question
addressed here.

\begin{rem}
The matching polynomial $M(G, x, y)$ of an arbitrary graph is not
polynomial reconstructible. The minimal counterexample for simple
graphs (with respect to the number of vertices and edges) are
constructed in \cite{LiShiTrinks}.
\end{rem}

The question arises here: whether or not there are such counterexamples
consisting of graphs with an arbitrary even number of vertices. In
\cite{LiShiTrinks}, we gave an affirmative answer to this question.

\section{Roots of beta-polynomials and independence polynomials}

Polynomials whose all zeros are real-valued numbers are said to be
{\it real}. Several graph polynomials have been known to be real; among
them the matching polynomial $\alpha(G,x)$ plays a distinguished
role \cite{godsilgutman,HeilmannLieb}.

Polynomials with only real roots arise in various applications in
control theory and computer science \cite{Visontai}, but also admit
interesting mathematical properties on their own. Newton noted that
the sequence of coefficients of such polynomials form a log-concave
(and hence unimodal) sequence. These polynomials also have strong
connections to totally positive matrices.

The fact that for all graphs, all zeros of the matching polynomial
are real-valued was first established by Heilmann and Lieb
\cite{HeilmannLieb}.

Let $C$ be a circuit contained in a graph $G$. If $C$ is a
Hamiltonian cycle, then $\alpha(G-C, x)\equiv 1$. In certain
considerations in theoretical chemistry
\cite{Aihara,LepovicGutmanPetrovicMizoguchi,Mizoguchi1,Mizoguchi2},
graph polynomials $\beta(G,C,x)$ are encountered, defined as
\begin{align}\label{eq2}
\beta(G,C,x)=\alpha(G,x)-2\alpha(G-C,x)
\end{align}
and
\begin{align}\label{eq3}
\beta(G,C,x)=\alpha(G,x)+2\alpha(G-C,x)
\end{align}
Formula \eqref{eq2} is used in the case of so-called H\"uckel-type
circuits, whereas formula \eqref{eq3} for the so-called M\"obius-type
circuits. For details, see \cite{Mizoguchi1}. These polynomials are
also called {\it circuit characteristic polynomials} \cite{Aihara}.

Already in the first paper devoted to this topic \cite{Aihara},
Aihara mentioned that the zeros of the $\beta$-polynomials are
real-valued, but gave no argument to support his claim. In the
meantime, for a number of classes of graphs it was shown that
$\beta(G,C, x)$ is indeed a real polynomial
\cite{Gutman3,Gutman4,GutmanMizoguchi,LepovicGutmanPetrovicMizoguchi,LiZhaoGutman,Mizoguchi2,LepovicGutmanPetrovic}.
In addition to this, by means of extensive computer searches not a
single graph with non-real $\beta$-polynomial could be detected. The
following conjecture has been put forward by Gutman and Mizoguchi in
\cite{Gutman3,Gutman4,GutmanMizoguchi}.

\begin{con} For any circuit $C$ contained in any graph $G$, the $\beta$-polynomials
 $\beta(G,C,x)$, Eqs. \eqref{eq2} and \eqref{eq3}, are real. \end{con}

Many results have been obtained. In particular, $\beta(G,C,x)$ has
been shown to be real for unicyclic graphs \cite{GutmanMizoguchi},
bicyclic graphs \cite{Mizoguchi2}, graphs in which no edge belongs
to more than one circuit \cite{Mizoguchi2}, graphs without
3-matchings (i.e., $m(G,3)=0$) \cite{LepovicGutmanPetrovicMizoguchi},
several (but not all) classes of graphs without 4-matchings
(i.e., $m(G,4)=0$) \cite{LepovicGutmanPetrovic}.

In \cite{LiGutmanMilovanovic}, Li et al. showed that the conjecture
is true for complete graphs. Actually, they proved a stronger result
for complete graphs.
\begin{thm}
For any circuit $C$ in the complete graph $K_n$, the polynomial $$\beta (K_n, C, t; x)=
\alpha (K_n, x)+t\alpha (K_n-C, x)$$ is real for any real $t$ such that
$|t|\leq n-1$.
\end{thm}
The proof offered in \cite{LiGutmanMilovanovic} relies on an earlier published theorem by
Tur\'an. In \cite{LiGutman2000}, Li and Gutman presented an elementary
self-contained proof for complete graphs. Finally, in
\cite{LiZhaoWang}, Li et al. showed that the conjecture is true for
all graphs, and therefore completely solved this conjecture.

\begin{thm}
For any circuit $C$ contained in any graph $G$, all roots of the
polynomial $\beta(G,C, x)$ are real.
\end{thm}

Chudnovsky and Seymour \cite{ChudnovskySeymour} proved the following result for
independence polynomial.
\begin{thm}\label{thm1}
If $G$ is clawfree, then all roots of its independence polynomial
are real.
\end{thm}

Theorem \ref{thm1} extends a theorem of \cite{HeilmannLieb},
answering a question posed by Hamidoune \cite{Hamidoune} and Stanley
\cite{Stanley}. Since all line graphs are clawfree, this extends the
result of \cite{HeilmannLieb}. Later, Levit and Mandrescu studied
the roots of independence polynomials of almost all very
well-covered graphs \cite{LevitMadrescu1}. In \cite{Mandrescu},
Mandrescu showed that starting from a graph $G$ whose independence
polynomial has only real roots, one can build an infinite family of
graphs, whose independence polynomials have only real roots.

Real roots of other graph polynomials have also been extensively studied, such as
edge-cover polynomial \cite{AkbariOboudi}, the expected independence
polynomial \cite{BrownDilcherManna}, domination polynomial
\cite{BrownTufts}, sigma-polynomial \cite{ZhaoLiZhangLiu},
chromatic polynomial \cite{DongKoh,Jackson,Woodall}, Wiener
polynomial \cite{DehmerIlic}, flow polynomial \cite{Jackson},
Tutte polynomial \cite{EM, Verg}, etc. For
more results on the roots of graph polynomials, we refer to
\cite{DehmerShiMowshowitz,Haglund,HaglundOnoWagner,MakowskyRavveBlanchard,Nijenhuis,SavageVisontai,Visontai}.

\section*{Acknowledgment}
The authors are supported by NSFC and the ``973" program.

\end{document}